 \documentclass[10pt,twoside]{article}
\newcount\VOL\VOL=3\newcount\YEAR\YEAR=1998
\def\firstpage{1}\def\lastpage{1000}
\def\received{}\def\revised{}
\def\communicated{}
\font\eightrm=cmr8
\font\caps=cmcsc10                    
\font\Caps=cmcsc10 scaled \magstep1   

\makeatletter
\@specialpagefalse\headheight=8.5pt
\def\DocMath{}
\renewcommand{\@evenhead}{%
    \ifnum\thepage>\lastpage\rlap{\thepage}\hfill%
    \else\rlap{\thepage}\slshape\leftmark\hfill\caps\SAuthor\hfill\fi}%
\renewcommand{\@oddhead}{%
    \ifnum\thepage=\firstpage{\DocMath\hfill\llap{\thepage}}%
    \else{\slshape\rightmark}\hfill\caps\STitle\hfill\llap{\thepage}\fi}%
\makeatother

\def\TSkip{\bigskip}
\newbox\TheTitle{\obeylines\gdef\GetTitle #1
\ShortTitle  #2
\SubTitle    #3
\Author      #4
\ShortAuthor #5
\EndTitle
{\setbox\TheTitle=\vbox{\baselineskip=20pt\let\par=\cr\obeylines%
\halign{\centerline{\Caps##}\cr\noalign{\medskip}\cr#1\cr}}%
        \copy\TheTitle\TSkip\TSkip%
\def\next{#2}\ifx\next\empty\gdef\STitle{#1}\else\gdef\STitle{#2}\fi%
\def\next{#3}\ifx\next\empty%
    \else\setbox\TheTitle=\vbox{\baselineskip=20pt\let\par=\cr\obeylines%
    \halign{\centerline{\caps##} #3\cr}}\copy\TheTitle\TSkip\TSkip\fi%
\centerline{\caps #4}\TSkip\TSkip%
\def\next{#5}\ifx\next\empty\gdef\SAuthor{#4}\else\gdef\SAuthor{#5}\fi%
\ifx\received\empty\relax
    \else\centerline{\eightrm Received: \received}\fi%
\ifx\revised\empty\TSkip%
    \else\centerline{\eightrm Revised: \revised}\TSkip\fi%
\ifx\communicated\empty\relax
    \else\centerline{\eightrm Communicated by \communicated}\fi\TSkip\TSkip%
\catcode'015=5}}\def\Title{\obeylines\GetTitle}
\def\Abstract{\begingroup\narrower
    \parskip=\medskipamount\parindent=0pt{\caps }}
\def\EndAbstract{\par\endgroup\TSkip\TSkip}
\newbox\TheAdd\def\Addresses{\vfill\copy\TheAdd\vfill
    \ifodd\number\lastpage\vfill\eject\phantom{.}\vfill\eject\fi}
{\obeylines\gdef\GetAddress #1
\Address #2 
\Address #3
\Address #4
\EndAddress
{\def\xs{6truecm}
\setbox0=\vtop{{\obeylines\hsize=\xs#1}}\def\next{#2}
\ifx\next\empty 
     \setbox\TheAdd=\hbox to\hsize{\hfill\copy0\hfill}
\else\setbox1=\vtop{{\obeylines\hsize=\xs#2}}\def\next{#3}
\ifx\next\empty 
     \setbox\TheAdd=\hbox to\hsize{\hfill\copy0\hfill\copy1\hfill}
\else\setbox2=\vtop{{\obeylines\hsize=\xs#3}}\def\next{#4}
\ifx\next\empty\ 
     \setbox\TheAdd=\vtop{\hbox to\hsize{\hfill\copy0\hfill\copy1\hfill}
                \vskip20pt\hbox to\hsize{\hfill\copy2\hfill}}
\else\setbox3=\vtop{{\obeylines\hsize=\xs#4}}
     \setbox\TheAdd=\vtop{\hbox to\hsize{\hfill\copy0\hfill\copy1\hfill}
                \vskip20pt\hbox to\hsize{\hfill\copy2\hfill\copy3\hfill}}
\fi\fi\fi\catcode'015=5}}\gdef\Address{\obeylines\GetAddress}

\begin{document} 

\Title The Work of R.E. Borcherds
\ShortTitle The Work of R.E. Borcherds
\SubTitle   
\Author Peter Goddard

\ShortAuthor
\EndTitle

\Abstract 
{\it Laudation delivered at the International Congress of Mathematicians in
Berlin following the award of the Fields Medal to Richard Borcherds.}
\EndAbstract
\centerline{18 August 1998}
\Address
Peter Goddard

St John's College
Cambridge CB2 1TP
\Address
\Address
\Address
\EndAddress
%
%

\input amssym.def
\input amssym.tex

\def\refname{Selected Papers of R.E. Borcherds}
\def\Rop{{\Bbb R}}
\def\Cop{{\Bbb C}}
\def\Zop{{\Bbb Z}}
\def\L{{\cal L}}
\def\half{{1\over 2}}
\parskip=6pt
\parindent=0pt
\section{Introduction}
Richard Borcherds has used the study of certain exceptional and exotic 
algebraic structures to motivate the introduction of important new
algebraic concepts: vertex algebas and generalized Kac-Moody algebras,
and he has demonstrated their power by using them to prove the
``moonshine conjectures'' of Conway and Norton about the Monster Group and
to find whole new families of automorphic forms. 

A central thread in his research has been a particular Lie algebra,
now known as the Fake Monster Lie algebra, which is, in a certain
sense, the simplest known example of a generalized Kac-Moody algebra
which is not finite-dimensional or affine (or a sum of such algebras). 
As the name might suggest, this algebra {\it appears} to have something
to do with the Monster group, {\it i.e.} the largest sporadic finite
simple group. 

The story starts with the observation that the Leech
lattice can be interpreted as the Dynkin diagram for a
Kac-Moody  algebra, $\L_\infty$. But $\L_\infty$ is difficult to handle; its
root multiplicities are not known explicitly. Borcherds showed how to enlarge
it to obtain the more amenable 
Fake Monster Lie algebra. In order to construct this algebra, Borcherds
introduced the concept of a vertex algebra, in the process establishing a
comprehensive algebraic approach to (two-dimensional) conformal field theory, a
subject of major importance in theoretical physics in the last thirty years.

To provide a general context for the Fake Monster Lie algebra, Borcherds
has developed the theory of generalized Kac-Moody algebras, proving, in
particular, generalizations of the Kac-Weyl character and denominator formulae.
The denominator formula for the Fake Monster Lie algebra motivated Borcherds to
construct a ``real'' Monster Lie algebra, which he used to prove the
moonshine conjectures. The results for the Fake Monster Lie algebra
also motivated Borcherds to explore the properties of the  denominator formula
for other generalized Kac-Moody algebras, obtaining remarkable product
expressions for modular functions, results on the moduli spaces of certain
complex surfaces and much else besides.

\section{The Leech Lattice and the Kac-Moody Algebra $\L_\infty$}

We start by recalling that a
finite-dimensional simple complex Lie algebra, $\L$, can be expressed in
terms of generators and relations as follows. There is a non-singular
invariant bilinear form $(,)$ on $\L$ which induces such a form on the 
$\hbox{rank} \L$ dimensional space spanned by the roots of $\L$. Suppose
$\{\alpha_i : 1\leq i\leq \hbox{rank} \L\}$ is a basis of simple roots for $\L$.
Then the numbers $a_{ij}=(\alpha_i,\alpha_j)$ have the following
properties:
\begin{eqnarray}
a_{ii}&>&0, \label{km1}\\ 
a_{ij}&=&a_{ji}, \label{km2} \\ 
a_{ij}&\leq&0\quad \hbox{ if } i\ne j, \label{km3}\\ 
2a_{ij}/a_{ii}&\in&\Zop.  \label{km4}
\end{eqnarray}
The symmetric matrix $A=(a_{ij})$ obtained in this way is positive
definite.

The algebra $\L$ can be reconstructed from the matrix $A$ by the
system of generators and relations used to define $\L_\infty$,
\begin{equation}\label{s1} 
[e_i,f_i]=h_i,\qquad [e_i,f_j]=0\quad\hbox{ for } i\ne j,
\end{equation}
\begin{equation}\label{s2}
[h_i,e_j]=a_{ij}e_j,\qquad [h_i,f_j]=-a_{ij}f_j,
\end{equation}
\begin{equation}\label{s3}
{\rm Ad}(e_i)^{n_{ij}}(e_j)={\rm Ad}(f_i)^{n_{ij}}(f_j)=0,
\quad\hbox{for } n_{ij} = 1 - 2a_{ij}/a_{ii}.
\end{equation}
These relations can be used to define a 
Lie algebra, $\L_A$, for any matrix $A$ satisfying the conditions
(\ref{km1}-\ref{km4}).
$\L_A$ is called a (symmetrizable) Kac-Moody algebra. If $A$ is positive
definite, $\L_A$ is semi-simple and, if $A$ is positive semi-definite, 
$\L_A$ is a sum of affine and finite-dimensional algebras.

Although Kac and Moody only explicitly considered the situation in which
the number of simple roots was finite, the theory of Kac-Moody algebras
applies to algebras which have a infinite number of
simple roots. Borcherds and others \cite{b0} showed how to construct such an
algebra with simple roots labelled by the points of the Leech lattice,
$\Lambda_L$. We can conveniently describe $\Lambda_L$ as a subset of the unique
even self-dual lattice, ${\rm II}_{25,1}$, in 26-dimensional Lorentzian space,
$\Rop^{25,1}$. ${\rm II}_{25,1}$ is the set of points whose coordinates are all
either integers or half odd integers and which have integral inner product with
the  vector
$(\half,\ldots,\half;\half)\in\Rop^{25,1}$, where the norm of
$x = (x_1,x_2,\ldots,x_{25};x_0)$ is
$x^2 = x_1^2+x_2^2 + \ldots + x_{25}^2 - x_0^2.$

The vector $\rho=(0,1,2,\ldots,24;70)\in{\rm II}_{25,1}$ has zero norm,
$\rho^2=0$;  the Leech lattice can be shown to be isomorphic to the set
$\{x\in{\rm II}_{25,1}:x\cdot\rho=-1\}$ modulo
displacements by $\rho$. We can take the representative points for the
Leech lattice to have norm 2 and so obtain an isometric correspondence between
$\Lambda_L$ and $\{r\in{\rm II}_{25,1}:r\cdot\rho=-1, r^2=2\}$.
Then, with each point $r$ of the Leech lattice, we can associate a 
reflection $x\mapsto \sigma_r(x)=x- (r\cdot x) r$ which is an
automorphism of 
${\rm II}_{25,1}$. Indeed these reflections $\sigma_r$ generate a Weyl
group, $W$, and the whole automorphism group of ${\rm II}_{25,1}$
is the semi-direct product of $W$ and the automorphism group of 
the affine Leech lattice, which is the Coxeter/Dynkin diagram of the Weyl group
$W$.
To this Dynkin diagram can be associated an infinite-dimensional Kac-Moody
algebra,
$\L_\infty$, generated by elements $\{e_r,f_r,
h_r:r\in\Lambda_L\}$ subject to the relations (\ref{s1}-\ref{s3}). Dividing by
the linear combinations of the $h_r$ which are in the centre reduces its rank to
26.

The point about Kac-Moody algebras is that they share many of the
properties enjoyed by semi-simple Lie algebras. In particular, we
can define a Weyl group, $W$, and for
suitable ({\it i.e.} lowest weight) representations, there is a
straightforward generalization of the Weyl character formula.  
For a representation with lowest weight $\lambda$, this
generalization, the Weyl-Kac character formula, states
\begin{equation}\label{wk}
\chi_\lambda= \left.\sum_{w\in
W}\det(w)w(e^{\rho+\lambda})\right/e^\rho\prod_{\alpha  >0}
\left(1-e^{\alpha}\right)^{m_\alpha},
\end{equation}
where $\rho$ is the Weyl vector, with $\rho\cdot r = -r^2/2$ for
all simple roots $r$, $m_\alpha$ is the multiplicity of the root
$\alpha$, the sum is over the elements $w$ of the Weyl group $W$, and
the product is  over positive roots $\alpha$, that is roots which can be
expressed as the sum of a subset of the simple roots with positive
integral coefficients. 

Considering even just the trivial representation, for
which $\lambda=0$ and $\chi_0=1$, yields a potentially
interesting relation from (\ref{wk}),
\begin{equation}\label{df}
\sum_{w\in
W}\det(w)w(e^{\rho})=e^\rho\prod_{\alpha  >0}
\left(1-e^{\alpha}\right)^{m_\alpha}.
\end{equation}
Kac showed that this denominator identity produces the
Macdonald identities in the affine case. Kac-Moody algebras, other than
the finite-dimensional and affine ones, would seem to offer the prospect
of new identities generalizing these but the problem is that in other
cases of Kac-Moody algebras, although the simple roots are known (as
for
$\L_\infty$), which effectively enables the sum over the Weyl group to 
be evaluated, the root multiplicities, $m_\alpha$,
are not known, so that the product over positive roots can not be evaluated. 

No general simple explicit formula is known for the root multiplicities
of $\L_\infty$ but, using the ``no-ghost'' theorem of string theory,
I. Frenkel established the bound $m_\alpha \leq p_{24}(1-\half
\alpha^2)$, where $p_k(n)$ is the number of partitions of $n$ using
$k$ colours.  This bound is saturated for some of the roots of 
$\L_\infty$ and, where it is not, there is the impression that that is
because something is missing. What seems to be missing are some simple roots
of zero or negative norm. In Kac-Moody algebras all the simple roots are 
specified by (\ref{km1}) to be of positive norm, even
though some of the other roots they generate may not be.

\section{Vertex Algebras}

Motivated by Frenkel's work, Borcherds introduced  \cite{b2} the definition of a
vertex algebra, which could in turn be used to define Lie algebras with root
multiplicities which are explicitly calculable. A vertex algebra is a 
graded complex vector space, $V=\bigoplus_{n\in\Zop}V_n$, together
with a ``vertex operator'', $a(z)$, for each $a\in V$, which is a formal power
series in the complex variable $z$, 
\begin{equation}
a(z) = \sum_{m\in\Zop} a_mz^{-m-n},\qquad\hbox{for } a\in V_n,
\end{equation}
where the operators $a_m$ map $V_n\rightarrow V_{n-m}$ and 
satisfy the following properties:
\begin{enumerate}
\item
$a_nb=0$ for $n>N$ for some integer $N$ dependent on $a$ and $b$;
\item
there is an operator (derivation) $D:V\rightarrow V$ such that
$[D,a(z)]={d\over dz}a(z)$;
\item
there is a vector ${\bf 1}\in V_0$ such that ${\bf 1}(z) = 1$, 
$D{\bf 1}=0$;
\item
$a(0){\bf 1} = a$;
\item
$(z-\zeta)^N\left(a(z)b(\zeta)-b(\zeta)a(z)\right)=0$
for some integer $N$ dependent on $a$ and $b$.\label{local}
\end{enumerate}
[We may define vertex operators over other fields or over the integers with
more effort if we wish but the essential features are brought out in the
complex case.]

The motivation for these axioms comes from string theory, where the vertex 
operators describe the interactions of ``strings'' (which are to be interpreted
as models for elementary particles). Condition (\ref{local}) states that
$a(z)$ and $b(\zeta)$ commute apart from a possible pole at $z=\zeta$, {\it
i.e.} they are local fields in the sense of quantum field theory. A key result
is that, in an appropriate sense,
\begin{equation}
(a(z-\zeta)b)(\zeta)=a(z)b(\zeta)=b(\zeta)a(z).
\end{equation}
More precisely
\begin{equation}
\int_0d\zeta\int_\zeta dz\; (a(z-\zeta)b)(\zeta)f=
\int_0 dz\int_0d\zeta\; a(z)b(\zeta)f-\int_0d\zeta\int_0 dz\;  b(\zeta)a(z)f.
\label{duality}\end{equation}
where $f$ is a polynomial in $z$, $\zeta$, $z-\zeta$ and their inverses, and the
integral over
$z$ is a circle about $\zeta$ in the first integral, one about $\zeta$ and the
origin in the second integral and a circle about the origin excluding the
$\zeta$ in the third integral. The axioms originally proposed by Borcherds
\cite{b2} were somewhat more complicated in form and follow from those given
here from the conditions generated by (\ref{duality}). 

We can associate a vertex algebra to any even lattice $\Lambda$, the space
$V$ then having the structure of the tensor product of the complex group ring 
$\Cop(\Lambda)$ with the symmetric algebra of a sum $\bigoplus_{n>0}\Lambda_n$
of copies $\Lambda_n, n\in\Zop,$ of $\Lambda$. In terms of string theory, this is
the Fock space describing  the (chiral) states of a string moving in a space-time
compactified into a torus by imposing perodicity under displacements by the
lattice
$\Lambda$.

The first triumph of vertex algebras was to provide a natural setting for the
Monster group, $M$. $M$ acts on a graded
infinite-dimensional space
$V^\natural$, constructed by Frenkel, Lepowsky and Meurman, where
$V^\natural=\bigoplus_{n\geq -1}V^\natural_n$, and the dimension of $\dim
V^\natural_n$ is the coefficent, $c(n)$ of $q^n$ in the elliptic modular
function,
\begin{equation}
j(\tau)-744 = \sum_{n=-1}^\infty c(n)q^n= q^{-1} +
196884q+21493760q^2+\ldots,\qquad q=e^{2\pi i \tau}.
\end{equation} 
A first thought might have been that the Monster group should be
related to the space $V_{\Lambda_L}$, the vertex algebra directly associated
with the Leech lattice, but $V_{\Lambda_L}$ has a grade $0$ piece of dimension 24
and the lowest non-trivial representation of the Monster is of dimension 196883. 
$V^\natural$ is related to $V_{\Lambda_L}$ but is a sort of twisted version of
it; in string theory terms it corresponds to the string moving on an orbifold
rather than a torus. 

The Monster group is precisely the group of automorphisms of the vertex algebra
$V^\natural$,
\begin{equation}
ga(z)g^{-1}=(ga)(z), \qquad g\in M.
\end{equation}
This characterizes $M$ in a way similar to the way that two other sporadic
simple finite groups, Conway's group $Co_1$ and the Mathieu group $M_{24}$,
can be characterized as the automorphism groups of the Leech lattice (modulo
$-1$) and the Golay Code, respectively.

\section{Generalized Kac-Moody Algebras}

In their famous moonshine conjectures, Conway and Norton went far
beyond  the existence of the graded representation $V^\natural$ with dimension
given by
$j$. Their main conjecture was that, for each element $g\in M$, the Thompson
series
\begin{equation}
T_g(q) = \sum_{n=-1}^\infty \hbox{Trace}(g|V^\natural_n)q^n
\end{equation} 
is a Hauptmodul for some genus zero subgroup, $G$, of $\hbox{SL}_2(\Rop)$,
{\it i.e.}, if $H=\{\tau:
\hbox{Im}(\tau)>0\}$ denotes the upper half complex plane, $G$ is such that the
closure of
$H/G$ is a compact Riemann surface, $\overline{H/G}$, of genus zero with a finite
number of points removed and 
$T_g(q)$ defines an isomorphism of $\overline{H/G}$ onto the Riemann sphere.  

To attack the moonshine conjectures it is necessary to introduce some Lie
algebraic structure. For any vertex algebra, $V$, we can introduce \cite{b2,b8} a
Lie algebra  of operators 
\begin{equation}
L(a) = {1\over 2\pi i}\oint a(z) dz = a_{-h+1}, \qquad a\in V_h,
\end{equation}
Closure $[L(a),L(b)]= L(L(a)b)$ follows from (\ref{duality}), but 
this does not define a Lie algebra structure directly on $V$ because 
$L(a)b$ is not itself antisymmetric in $a$ and $b$. However, $DV$ is in the
kernel of the map $a\mapsto L(a)$ and $L(a)b=-L(b)a$ in $V/DV$, so it does define
a Lie algebra $\L^0(V)$  on this quotient \cite{b2}, but this is not
the  most interesting Lie algebra associated with $V$.

Vertex algebras of interest come with an additonal structure, an action of the
Virasoro algebra, a central extension of the Lie algebra of polynomial vector 
fields on the circle, spanned by $L_n, n\in\Zop$ and $1$,
\begin{equation}
[L_m,L_n]=(m-n)L_{m+n} + {c\over 12}m(m^2-1)\delta_{m,-n},\qquad [L_n,c]=0,
\end{equation}
with $L_{-1}=D$ and $L_0a=ha$ for $a\in V_h$. For $V_\Lambda$, $c=\dim \Lambda$,
and for $V^\natural$, $c=24$. The Virasoro algebra plays a central 
role in string theory. The space of ``physical states'' of the string is
defined by the Virasoro conditions: let
$P^k(V)=\{a\in V : L_0a=ka; L_na=0, n>0\}$, the space of physical states is
$P^1(V)$. The space $P^1(V)/L_{-1}P^0(V)$ has a Lie algebra structure defined
on it (because $L_{-1}V\cap P^1(V)\subset L_{-1}P^0(V)$). This can be reduced 
in size further using a contravariant form (which it possesses naturally for 
lattice theories). The ``no-ghost'' theorem states that the
space of physical states $P^1(V)$ has lots of null states and is positive 
semi-definite for $V_\Lambda$, where $\Lambda$ is a Lorentzian lattice with
$\dim\Lambda\leq 26$. So we can quotient $P^1(V)/L_{-1}P^{0}(V)$ further by 
its null space with the respect to the contravariant form to obtain a 
Lie algebra $\L(V)$. 

The results of factoring by the null space are most dramatic when $c=26$. The
vertex algebra $V_L$ has a natural grading by the lattice $L$ and the
``no-ghost'' theorem states that the dimension of the subspace of 
$\L(V)$ of non-zero grade $\alpha$ is $p_{24}(1-\half\alpha^2)$ if $\Lambda$ is
a Lorentzian lattice of dimension 26 but
$p_{k-1}(1-\alpha^2/2)-p_{k-1}(\alpha^2/2)$ if $\dim\Lambda=k\ne 26$, $k>2$.
Thus the algebra $\L_M'=\L(V_{{\rm II}_{25,1}})$ saturates Frenkel's bound, and
Borcherds initially named it the ``Monster Lie algebra'' because it appeared to
be directly connected to the Monster; it is now known as the ``Fake Monster Lie
algebra''. 

Borcherds \cite{b8} had the great insight not only to construct the Fake Monster
Lie algebra, but also to see how to generalize the definition of a Kac-Moody
algebra effectively in order to bring $\L_M'$ within the fold. What was required
was to  relax the condition (\ref{km1}), requiring roots to have positive norm,
and to allow them to be either zero or negative norm.
The condition (\ref{km4}) then needs modification to apply only in the
space-like case $a_{ii}>0$ and the same applies to the condition (\ref{s3}) on
the generators. The only condition which needs to be added is  that
\begin{equation}
[e_i,e_j]=[f_i,f_j]=0\qquad\hbox{if } a_{ij}=0.
\end{equation}

The closeness of these conditions to those for Kac-Moody algebras means that
most of the important structural results carry over; in particular there is
a generalization of the Weyl-Kac character formula for
representations with highest weight $\lambda$,

\begin{equation}\label{gwk}
\chi_\lambda= \left.\sum_{w\in
W}\det(w)w\left(e^\rho\sum_\mu  \epsilon_\lambda(\mu)e^{\mu+\lambda}\right)
\right/e^\rho\prod_{\alpha >0}
\left(1-e^{\alpha}\right)^{m_\alpha},
\end{equation}
where the second sum in the numerator is over vectors $\mu$ and
$\epsilon_\lambda(\mu)=(-1)^n$ if $\mu$ can be expressed as ths sum of 
$n$ pairwise orthogonal simple roots with non-positive norm,
all orthogonal to $\lambda$, and $0$ otherwise. Of course, putting $\lambda=0$
and 
$\chi_\lambda=1$ again gives a denominator formula.

The description of generalized Kac-Moody algebras in terms of generators and
relations enables the theory to be taken over rather simply from that of
Kac-Moody algebras but it is not so convenient as a method of recognising them 
in practice, {\it e.g.} from amongst the algebras $\L(V)$ previously constructed
by Borcherds. But Borcherds \cite{b4} gave an alternative characterization of
them as as graded algebras with an ``almost postitive definite'' contravariant
bilinear form. More precisely, he showed that a graded Lie algebra,
$\L=\bigoplus_{n\in\Zop}\L_n$, is a generalized Kac-Moody algebra if 
the following conditions are satisfied:
\begin{enumerate}
\item
$\L_0$ is abelian and $\dim \L_n$ is finite if $n\ne0$;
\item
$\L$ possesses an invariant bilinear form such that $(\L_m,\L_n)=0$ if $m\ne n$;
\item
$\L$ possesses an involution $\omega$ which is $-1$ on $\L_0$ and such that 
$\omega(\L_m)\subset \L_{-m}$;
\item
the contravariant bilinear form $\langle L, M\rangle=-(L,\omega(M))$  is
positive definite on $\L_m$ for $m\ne 0$ ;
\item
$\L_0\subset [\L,\L]$.
\end{enumerate}

This characterization shows that the Fake Monster Lie algebra, $\L_M'$, is a
generalised Kac-Moody algebra, and its root multiplicities are known to be
given by $p_{24}(1-\half\alpha^2)$, but Borcherds' theorem establishing the
equivalence of his two definitions does not give a constructive method of finding
the simple roots. As we remarked in the context of Kac-Moody algebras, if we knew
both the root multiplicities and the simple roots, the denominator formula 
\begin{equation}\label{gdf}
\sum_{w\in
W}\det(w)w\left(e^\rho\sum_\mu  \epsilon_\mu(\alpha)e^{\mu}\right)
=e^\rho\prod_{\alpha 
>0}
\left(1-e^{\alpha}\right)^{m_\alpha}
\end{equation}
might provide an interesting identity. Borcherds solved \cite{b8} the problem of
finding the simple roots, or rather proving that the obvious ones were all that
there were, by inverting this argument. The positive norm simple roots can be
identified with the Leech lattice as for $\L_\infty$. Writing
${\rm II}_{25,1}=\Lambda_L\oplus{\rm II}_{1,1}$, which follows by uniqueness or
the earlier comments, the `real' or space-like simple roots are 
$\{(\lambda, 1, \half\lambda^2-1):\lambda\in\Lambda_L\}$. (Here we are using we
are writing
${\rm II}_{1,1}=\{(m,n): m,n\in\Zop\}$ with $(m,n)$ having norm $-2mn$.)
Light-like simple roots are quite easily seen to be $n\rho$, where $n$ is a
positive integer and $\rho=(0,0,1)$. The denominator identity is then used to
prove that there are no other light-like and that there are no time-like simple
roots. 

The denominator identity provides a remarkable relation between modular
functions (apparently already known to some of the experts in the
subject) which is the precursor of other even more remarkable identities. If we
restrict attention to vectors
$(0,\sigma,\tau)\in
{\rm II}_{25,1}\otimes\Cop$, with $\hbox{Im}(\sigma)>0$, $\hbox{Im}(\tau)>0$, it
reads
\begin{equation}\label{fmid}
p^{-1}\prod_{m>0,n\in\Zop}(1-p^mq^n)^{c'(mn)}
=\Delta(\sigma)\Delta(\tau)(j(\sigma)-j(\tau))
\end{equation}
where $c'(0)=24$, $c'(n) = c(n)$ if $n\ne 0$, $p= e^{2\pi i\sigma}$,
$q=e^{2\pi i\tau}$, and
\begin{equation}
\Delta(\tau)^{-1}= q^{-1}\prod_{n>1}(1-q^n)^{-24}=\sum_{n\geq 0}p_{24}(n)q^{n-1}.
\end{equation}

\section{Moonshine, the Monster Lie Algebra and Automorphic Forms}

The presence of $j(\sigma)$ in (\ref{fmid}) suggests a relationship to the 
moonshine conjectures and Borcherds used \cite{b10, b11} this as motivation to
construct the ``real'' Monster Lie Algebra, $\L_M$ as one with denominator
identity obtained by multiplying each side of (\ref{fmid}) by
$\Delta(\sigma)\Delta(\tau)$, to obtain the simpler formula
\begin{equation}\label{mid}
p^{-1}\prod_{m>0,n\in\Zop}(1-p^mq^n)^{c(mn)}
=j(\sigma)-j(\tau).
\end{equation}
This looks like the denominator formula for a generalised Kac-Moody algebra 
which is graded by ${\rm II}_{1,1}$ and is such that the dimension of the
subspace of grade $(m,n)\ne (0,0)$ is $c(mn)$, the dimension of
$V^\natural_{mn}$. It is not difficult to see that this can be constructed by
using the vertex algebra which is the tensor product $V^\natural\otimes V_{{\rm
II}_{1,1}}$ and defining $\L_M$ to be the generalised Lie algebra,
$\L(V^\natural\otimes V_{{\rm II}_{1,1}})$, constructed from the physical states.

Borcherds used \cite{b10, b11} twisted forms of the denominator identity for
$\L_M$ to prove the moonshine conjectures. The action of $M$ on $V^\natural$
provides an action on $V=V^\natural\otimes V_{{\rm II}_{1,1}}$ induces an action
on the physical state space $P^1(V)$ and on its quotient, $\L_M=\L(V)$, by its
null space. The ``no-ghost'' theorem implies that the part of $\L_M$ of grade
$(m,n)$, $(\L_M)_{(m,n)}$, is isomorphic to
$V^\natural_{mn}$ as an $M$ module. Borcherds adapted the  argument he used to
establish the denominator identity to prove the twisted relation
\begin{eqnarray}
p^{-1}&\hskip-42pt\exp\left(-\sum_{N>0}\sum_{m>0,n\in\Zop}{\rm Tr}
(g^N|V^\natural_{mn})p^{mN}q^{nN}/N\right) \nonumber\\
&\hskip60pt=\sum_{m\in\Zop}{\rm Tr}(g|V^\natural_{m})p^m
-\sum_{n\in\Zop}{\rm Tr}(g|V^\natural_{n})q^n.
\end{eqnarray}
These relations on the Thompson series are sufficient to determine them from
their first few terms and to establish that they are modular functions of genus
0.

Returning to the Fake Monster Lie Algebra, the denominator formula given
in (\ref{fmid}) was restricted to vectors of the form $v=(0,\sigma,\tau)$ but
we consider it for more general $v\in{\rm II}_{25,1}\otimes\Cop$, giving the
denominator function
\begin{equation}
\Phi(v)=\sum_{w\in W}\det(w)e^{2\pi i(w(\rho),v)}\prod_{n>0}
\left(1-e^{2\pi in(w(\rho),v)}\right)^{24}.
\end{equation}  
This expression converges for $\hbox{Im}(v)$ inside a certain cone (the positive
light cone). Using the explicit form for
$\Phi(v)$ when
$v=(0,\sigma,\tau)$, the known properties of $j$ and $\Delta$ and the fact that
$\Phi(v)$ manifestly satisfies the wave equation, Borcherds
\cite{b11,b13,b18} establishes that 
$\Phi(v)$ satisfies the functional equation
\begin{equation}
\Phi(2v/(v,v))=-((v,v)/2)^{12}\Phi(v).
\end{equation}
It also has the properties that $\Phi(v+\lambda)=\Phi(v)$ for
$\lambda\in{\rm II}_{25,1}$ and $\Phi(w(v))=\det(w)\Phi(v)$ for $w\in
\rm{Aut}({\rm II}_{25,1})^+$, the group of automorphisms of the lattice
${\rm II}_{25,1}$ which preserve the time direction. These transformations
generate a discrete subgroup of the group of conformal transformations on
$\Rop^{25,1}$, which is itself isomorphic to  $O_{26,2}(\Rop)$; in fact the
discrete group is isomorphic to $\rm{Aut}({\rm II}_{26,2})^+$. The denominator
function for the Fake Monster Lie algebra defines in this way an automorphic
form of weight 12 for the discrete subgroup $\rm{Aut}({\rm II}_{26,2})^+$ of 
 $O_{26,2}(\Rop)^+$. This result once obtained is seen not to depend essentially
on the dimension 26 and Borcherds has developed this approach of obtaining
representations of modular functions as infinite products from denominator
formulae for generalized Kac-Moody algebras to obtain a
plethora of beautiful formulae \cite{b13,b18,b20}, {\it e.g.}
\begin{equation}
j(\tau)=q^{-1}\prod_{n>0}(1-q^n)^{c_0(n^2)}=q^{-1}(1-q)^{-744}(1-q^2)^{80256}
(1-q^3)^{-12288744}\ldots,
\end{equation}
where $f_0(\tau)=\sum_nc_0(n)q^n$ is the unique modular form of weight $\half$
for  for the group $\Gamma_0(4)$ which is such
that
$f_0(\tau) = 3q^{-3}+{\cal O}(q)$ at
$q=0$ and $c_0(n) = 0$ if $n\equiv 2$ or 3 mod 4.
He has also used
these denominator functions to establish results about the moduli spaces of
Enriques surfaces and and families of K3 surfaces \cite{b17,b19}.

Displaying penetrating insight, formidable technique and brilliant originality,
Richard Borcherds has used the beautiful properties of some exceptional
structures to motiviate new algebraic theories of great power with profound
connections with other areas of mathematics and physics. He has used them to
establish outstanding conjectures and to find new deep results in classical
areas of mathematics. This is surely just the beginning of what we have to learn
from what he has created.

\Addresses

\begin{thebibliography}{10}


\bibitem{b0} {\it A monster Lie algebra?} (with J.H. Conway, L. Queen and
N.J.A. Sloane)  Adv. Math. {\bf 53} (1984) 75-79.  


\bibitem{b2} {\it Vertex algebras, Kac-Moody algebras and the monster,} 
Proc. Nat. Acad. Sci. U.S.A. {\bf 83} (1986) 3068-3071. 



\bibitem{b4} {\it Generalized Kac-Moody algebras,}
J. Alg. {\bf 115} (1988) 501-512. 




\bibitem{b8} {\it The monster Lie algebra,}
Adv. Math. {\bf 83} (1990) 30-47. 


\bibitem{b10} {\it Monstrous moonshine and monstrous Lie algebras,}
Invent. Math. {\bf 109} (1992) 405-444.

\bibitem{b11} {\it Sporadic groups and string theory,} in
{\it Proceedings of the
First European Congress of
Mathematics,  Paris July 1992}, ed. A. Joseph {\it et al.}, Vol. 1,
 Birkhauser (1994) pp. 411-421.

\bibitem{b13} {\it Automorphic forms on $O_{s+2,2}(R)$ and infinite
products,}  Invent. Math. {\bf 120} (1995) 161-213.


\bibitem{b17} {\it The moduli space of Enriques surfaces and the 
fake monster Lie superalgebra,} Topology {\bf 35} (1996) 699-710. 

\bibitem{b18} {\it Automorphic forms and Lie algebras,} in Current
developments in mathematics, International Press (1996).


\bibitem{b19} {\it Families of K3 surfaces,} (with L. Katzarkov, T.
Pantev, N.I.  Shepherd-Barron) J. Algebraic Geometry {\bf 7} (1998)
183-193. 

\bibitem{b20} {\it Automorphic forms with singularities on
Grassmannians,} Invent. Math. {\bf 132} (1998) 491-562.

\end{thebibliography}
\end{document}